\documentclass{article}
\usepackage{amssymb,amsmath,amsthm}

\begin{document}

\large

\theoremstyle{plain}
\newtheorem{theorem}{Theorem}
\newtheorem{lemma}[theorem]{Lemma}
\newtheorem{proposition}[theorem]{Proposition}
\newtheorem{corollary}[theorem]{Corollary}
\newtheorem*{main}{Main Theorem}

\theoremstyle{definition}
\newtheorem*{definition}{Definition}

\theoremstyle{remark}
\newtheorem*{example}{Example}
\newtheorem*{remark}{Remark}
\newtheorem*{remarks}{Remarks}

\newcommand{\sse}{\subseteq}

\newcommand{\qq}{{\mathbb Q}}
\newcommand{\rr}{{\mathbb R}}
\newcommand{\nn}{{\mathbb N}}
\newcommand{\zz}{{\mathbb Z}}

\newcommand{\al}{\alpha}
\newcommand{\be}{\beta}
\newcommand{\ep}{\epsilon}
\newcommand{\la}{\lambda}
\newcommand{\de}{\delta}

\title{Lattice polytopes with distinct pair-sums}

\author{M. D. Choi \\ University of Toronto \\Toronto, M5S 3G3, Canada\\
\\ T. Y. Lam \\ University of California\\Berkeley, CA 94720 \\
\\Bruce Reznick  \\ University of Illinois\\ Urbana, IL 61801}

\maketitle

Let $\cal P$ be a lattice polytope in $\rr^n$, the convex
hull of a finite set in $\zz^n$, and let
$$
{\cal L}({\cal P}): = {\cal P} \cap \zz^n = \{v_1,\dots,v_N\}, 
$$
where $N = N({\cal P}): = |{\cal L}({\cal P}) |$. 
Suppose the $N +
\binom N2$ points in ${\cal L}({\cal P})+ {\cal L}({\cal P})$,
$$
2v_1, \dots,2v_N;\ v_1+v_2,\ v_1+v_3, \dots, v_{N-1}+v_N
$$
are distinct. In this case, we say that ${\cal P}$ is a {\it distinct
pair-sum\/} or {\it dps\/} polytope. Our interest in dps polytopes
comes from the study of the representation of polynomials as a sum of
squares of polynomials.

The following lemma offers two other geometrical characterizations of 
dps polytopes.
\begin{lemma}
Let $\cal P$ be a lattice polytope. Then the following are equivalent: 

\smallskip\noindent
$(1)$ ${\cal L}({\cal P})$ is a dps polytope. \\
$(2)$ ${\cal L}({\cal P})$ does not contain the vertices of a
(nondegenerate) parallelogram, and does not contain three collinear
points.\\
$(3)$ Suppose $\,v\neq v'\,$ and $\,w\neq w'\,$ are in ${\cal L}({\cal P})$.
Then $\,v'-v\,$ and $\,w'-w\,$ are parallel only if 
$\,\{v,v'\}=\{w,w'\}$.
\end{lemma}
\begin{proof} $(1)\Rightarrow (2)$.
Suppose $v_1, v_2, v_3, v_4 \in{\cal L}({\cal P})$ are the vertices of
a parallelogram. Then $v_1 - v_2 = v_3 - v_4$ implies $v_1 + v_4 = v_2
+ v_3$, so that ${\cal P}$ is not dps. Now suppose $v_1, v_2, v_3 \in{\cal
L}({\cal P})$, and $v_2$ is interior to the line segment
$\overline{v_1v_3}$. If $v_2$ is the midpoint of the segment, then
$v_2 + v_2 = v_1 + v_3$, so ${\cal P}$ is not dps. Otherwise, we may
assume that $v_2$ is closer to $v_1$ than to $v_3$. Then $v_4 = v_2
+(v_2 - v_1)$ will also be a lattice point on the 
line segment $\overline{v_1v_3}$, and $v_2$ is the midpoint of
$\overline{v_1v_4}$; again, 
${\cal P}$ is not dps. 

\smallskip\noindent
$(2)\Rightarrow (3)$. For $\,u\in {\mathbb Z}^n$, let $\,g(u)
=\mbox{gcd}(u_1,\dots,u_n)$. Suppose $\,g(u'-u)=d>1$. Then
$\,u'-u=du''\,$ for $\,u''\in {\mathbb Z}^n$, and the line segment
$\,\overline{u\,u'}\,$ contains the lattice points $\,u,\,u+u'',\,\dots,
\,u+du''=u'$. Thus, if (2) holds and $\,u,u'\in {\cal L}({\cal P})$,
$\,u\neq u'$, we have $\,g(u'-u)=1$.  Suppose $\,w'-w=\alpha\cdot (v'-v)$.
Then $\,\alpha=p/q\,$ for nonzero integers $\,p,\,q$, and $\,q(w'-w)=
p(v'-v)$.  Hence $\,|q|=g(q(w'-w))=g(p(v'-v))=|p|$, so $\,\alpha =\pm 1$.
Now the parallelogram condition in (2) implies that
$\,\{v,v'\}=\{w,w'\}$.

\smallskip\noindent
$(3)\Rightarrow (1)$. If (3) holds for $\,{\cal P}$, and 
$\,v_i,v_j,v_k,v_{\ell}\in  {\cal L}({\cal P})\,$ with $\,i\notin
\{k,\ell\}$, then $\,v_i-v_k\neq  v_{\ell}-v_j$, and so
$\,v_i+v_{j}\neq v_k+v_{\ell}\,$.  This proves (1).
\end{proof}

Our main results are these: if ${\cal P}$ in $\rr^n$ is a dps polytope, then
$N({\cal P}) \le 2^n$, and, for every $n$, we 
construct dps polytopes in $\rr^n$ for which $N({\cal P}) = 2^n$. 

\medskip
\noindent \textbf{Example 1.}
Let ${\cal P} \subset \rr^2$ be the triangle with vertices
$\{(0,1),(1,2),(2,0)\}$. Then ${\cal P}$ is a dps polytope, because 
$$
{\cal L}({\cal P})= \{(0,1),\,(1,2),\,(2,0),\, (1,1)\},
$$
and 
$$
{\cal L}({\cal P})+ {\cal L}({\cal P}) =
\{(0,2),(1,2),(1,3),(2,1),(2,2),(2,3),(2,4),(3,1),(3,2),(4,0)\}. 
$$
We can view ${\cal P}$ as the
projection onto the first two coordinates of the triangle with
vertices $\{(0,1,2),(1,2,0),(2,0,1)\}$, 
which lies in the hyperplane $x_1+x_2+x_3 = 3$.
(In this example, we could have just as well taken the triangle with
vertices $\{(0,0),(1,2),(2,1)\}$; again, ${\cal L}({\cal P})$  will
consist of the vertices of $\,{\cal P}\,$ and $(1,1)$.)

\medskip
\noindent \textbf{Example 2.}
Let
$$
\displaylines{
{\cal A} = \{ (4,1,0,0), (0,4,1,0), (0,0,4,1), (1,0,0,4)\}, \cr
{\cal B} = \{(2,1,1,1), (1,2,1,1),
(1,1,2,1), (1,1,1,2)\};\cr}
$$
and let ${\cal P} = cvx({\cal A}\cup {\cal B})\subset \rr^4$ be the convex
hull of ${\cal A}\cup {\cal B}$. 
By construction, ${\cal P}$ is cyclically symmetric with respect to
its coordinates. 
It is not hard to show that  ${\cal L}({\cal P})= {\cal A}\cup {\cal B}$.
Suppose $w =
(w^{(1)},w^{(2)},w^{(3)},w^{(4)}) \in{\cal L}({\cal P}) $. Since $w$
is a convex combination of ${\cal A}\cup {\cal B}$, we have $w^{(i)}
\ge 0$ and $\sum_i w^{(i)}  = 5$. If $w^{(i)} \ge 1$ for all $i$, then $w$ 
must be a permutation of $(2,1,1,1)$ and so lies in $\cal
B$. Otherwise, $w^{(i)}  = 0$ for some $i$, and by cycling the coordinates,
we may assume that $w^{(1)} = 0$. But then
$w$ must be a convex combination of $(0,4,1,0)$ and $(0,0,4,1)$ and so
$w \in \cal A$.  A routine check, which we omit, shows that
the $8 + \binom 82 = 36$ sums in  ${\cal L}({\cal P}) +{\cal L}({\cal P})$ are
distinct. By projecting ${\cal P}$ onto its first three coordinates,
we obtain a dps polytope in $\rr^3$ with $N({\cal P}) = 8$.
\smallskip
\begin{theorem}
Suppose $\cal P$ is a dps polytope in $\rr^n$. Then $N({\cal P}) \le
2^n$.
\end{theorem}
\begin{proof}
If $N({\cal P}) > 2^n$, then by the Pigeonhole Principle, there exist
$v_i \neq v_j$ so that $v_i$ and $v_j$ are component-wise congruent
modulo 2. This means that $v_k = \frac 12(v_i + v_j)= v_i + \frac
12(v_j-v_i)$ is also a lattice point, and it follows from Lemma 1 that
$\cal P$ is not a dps polytope. 
\end{proof}

This argument is essentially the same one used to solve Putnam Problem
1971-A1 (see \cite{AKL}): ``Let there be given nine lattice points (points with
integral coordinates) in three dimensional Euclidean space. Show that
there is a lattice point on the interior of one of the line segments
joining two of these points."  The proof of Theorem 2 also applies to the
less restrictive class of convex polytopes which do not contain three
lattice points on a line. One such polytope is the $n$-cube ${\cal
C}_n = \{0,1\}^n$, which has many lattice parallelograms.

We shall say that a dps polytope ${\cal P} \subset \rr^n$ for which $N({\cal
P}) = 2^n$ is  {\it maximal}. The proof of Theorem 2 implies that no two
points in a dps polytope are component-wise congruent modulo 2; hence
a maximal dps polytope contains one representative from every
congruence class modulo 2 (and at most one representative from every
congruence class modulo $m$, $m \ge 3$).

Suppose $M$ is an $n \times n$ unimodular matrix with integer entries. Then $M$
defines a linear mapping on  $\rr^n$ (viewed as column vectors) by
matrix multiplication. Since linear mappings preserve inclusions and
both $M$ and 
$M^{-1}$ have integer entries, it is easy to see that ${\cal L}(M({\cal P}))
= M({\cal L}({\cal P}))$ for any lattice polytope ${\cal P}$, and
since linear mappings preserve sums, it is then clear that
${\cal P}$ is dps if and only if $M({\cal P})$ is dps. 

\begin{theorem}
 There exist maximal dps polytopes in $\rr^n$ for every $n$.
\end{theorem}
\begin{proof}
For $n=1$, let ${\cal P} = [0,1]$; for  $n =2,3$, consider Examples 1
and 2. Suppose now that  $\cal
P$ is a maximal dps polytope in $\rr^n$, $n \ge 3$. Write
$ {\cal L} = {\cal L}({\cal P})$ and define the (finite) set of differences
$$
{\cal D} = ({\cal L} - {\cal L})^*: = \{ v - v':\, v, v' \in {\cal L}, \,v
\neq v'\}. 
$$
Let $M$ be a unimodular integer matrix such that if $u \in {\cal D}$,
then $M(u) \notin {\cal D}$. (We shall construct such an $M$ below.)

We define the
polytope ${\cal P}'$ in $\rr^{n+1}$ as follows. Let 
$$
{\cal A} = \{ (v,0) \in \rr^{n+1}: v \in {\cal L}({\cal P})\},\qquad
{\cal B} = \{ (M(v),1) \in \rr^{n+1}: v \in {\cal L}({\cal P})\},
$$
and let ${\cal P}' = cvx({\cal A} \cup {\cal B})$.  
If $w = (w^{(1)},\dots,w^{(n+1)}) \in {\cal L}({\cal P}')$, then $0 \le
w^{(n+1)} \le 1$, hence $w^{(n+1)}$ equals 0 or 1. Thus, $w$ lies
either on the face determined by ${\cal A}$, in which case $w =
(v,0)$, or on the face determined by   ${\cal B}$, in which case $w =
(M(v),1)$. It follows that ${\cal L}({\cal P}') = {\cal A} \cup {\cal
B}$, so $N({\cal P}') = 2^{n+1}$.

Now consider $ {\cal L}({\cal P}') + {\cal L}({\cal P}')$; this
consists of three disjoint sets of points: 
$$\{(v_i,0) + (v_j,0)\}, \qquad \{(v_i,0) +
(M(v_j),1)\}, \qquad \{(M(v_i),1) + (M(v_j),1)\},
$$
where $v_i, v_j \in {\cal L}({\cal P})$. 
Since both ${\cal P}$ and $M({\cal P})$ are dps, the
sums in the first and the third set are distinct. 
For the second set, we suppose that
\begin{equation} \label{eq}
(v_i,0) + (M(v_j),1) = (v_k,0) + (M(v_\ell),1),
\end{equation}
or equivalently,
$$
v_i - v_k = M(v_\ell) - M(v_j) = M(v_\ell - v_j).
$$
If $j = \ell$, then $v_i - v_k = 0$, so $i = k$, which is the
only possible way for \eqref{eq} to hold in a dps polytope.
Otherwise, $j \neq \ell$, so  $M(v_\ell - v_j) = v_i - v_k \in {\cal D}$, a
contradiction to the 
choice of $M$. Thus, ${\cal P}'$ is a maximal dps polytope in $\rr^{n+1}$.

We now construct a matrix $M$ with the desired properties. First, let 
$$
R = \max\,\{ |u_j^{(k)}|: u_j \in {\cal D},\;1\leq k\leq n \}.
$$
and let $M$ be the $n \times n$ matrix given below:
$$
M =  \left( \begin{array}{ccccccc}
1+(R+1)^2 & R+1 & 0 & 0 & \dots & 0 & 0 \\
R+1 & 1 & R+1 & 0 & \dots & 0 & 0\\
0 & 0 & 1 & R+1 & \dots & 0 &0\\
0 & 0 & 0 & 1 & \dots & 0 &0\\
\dots & \dots & \dots &\dots & \dots & \dots &\dots\\
0 & 0 & 0 & 0 &\dots & 1 & R+1 \\
0 & 0 & 0 & 0 &\dots & 0 & 1
\end{array}
\right)
.
$$
(In words, the only non-zero entries in $M$ are the diagonal, the
superdiagonal, and the first entry in the second row.) It is easy to
see that $M$ is unimodular.  

We show now that for every $u \in \cal D$, at least one entry of
$w = M(u)$ has absolute value greater than $R$. This implies that $M(u)
\notin \cal D$, and will complete the proof. Write $u =
(u^{(1)},\dots,u^{(n)})$ and suppose that $k$ is the smallest index
such that $u^{(k)} \neq 0$. (Such an index exists because $0 \notin
\cal D$.)

If $k = 1$, then  $w^{(1)} = (1+(R+1)^2)u^{(1)} + (R+1)u^{(2)}$,
and hence
$$
|w^{(1)}| \ge |(1+(R+1)^2)u^{(1)}| -  (R+1)|u^{(2)}| \ge 1+(R+1)^2 -
R(R+1) = R+2.
$$
If $k\ge 2$,  then  $u^{(1)} = \dots = u^{(k-1)} = 0$, so
$w^{(k-1)} = (R+1)u^{(k)}$ and $|w^{(k-1)}| \ge R+1$.
Finally, we remark that  the same proof applies in the case $n=2$, if we take
as our matrix the $2 \times 2$ submatrix at the upper left of $M$.
\end{proof}

\noindent \textbf{Example 3.}
We illustrate the last construction by applying it to the polytope in
Example 1, for which
$$
{\cal D} = \{ \pm (0,1),\pm(1,-2), \pm (1,-1),\pm(1,0), \pm(1,1), \pm (2,-1)\},
$$
so $R = 2$ and 
$$
M = \left( \begin{array}{cc} 10 & 3 \\ 3 & 1 \end{array} \right).
$$
Thus, $cvx({\cal A} \cup {\cal B})$ is a maximal dps polytope in
$\rr^3$, where
$$
\displaylines{
{\cal A} = \{ (0,1,0), (1,1,0), (1,2,0), (2,0,0)\}, \cr
{\cal B} = \{(3,1,1), (13,4,1), (16,5,1), (20,6,1)\}.\cr}
$$
We could now apply the shear $(x_1,x_2,x_3) \mapsto
(x_1-3x_2-5x_3+5,x_2-x_3,x_3)$, which 
maps  ${\cal A}$ and ${\cal B}$ to 
$$
{\cal A}':=\  \{(2,1,0),(3,1,0),(0,2,0),(7,0,0)\}
$$
and
$$
{\cal B}':=\  \{(0,0,1),(1,3,1),(1,4,1),(2,5,1)\},
$$
respectively, in order to reduce the magnitude of the coordinates in
the example. 
\smallskip

Since any translate of a dps
polytope is also dps, we may always assume, as we have done in the
examples, 
that $\cal P$ lies in the non-negative orthant of $\rr^n$. In this
case, we define $s({\cal P})$, the {\it size} of $\cal P$:
$$
s({\cal P}) = \max \{ v_j^{(1)} + \cdots +  v_j^{(n)}: \ v_j \in {\cal
L}({\cal P})\}. 
$$
If $s = s({\cal P})$, then $\cal P$ can be viewed as a projection
onto the first $n$ coordinates of a polytope in $\rr^{n+1}$ which lies
in the simplex 
$$
\Delta_{n+1}(s):=\ \{ u = (u^{(1)},\dots, u^{(n+1)}): \  u^{(i)} \ge 0,
\ \sum_{i=1}^{n+1} u^{(i)} = s\}. 
$$

Let $s_n$ denote the
minimum size of any maximal dps polytope in $\rr^n$. Examples 1 and 2
show that   $s_2 \le 3$ and $s_3 \le 5$. It is not difficult
to show that these estimates are sharp. The first
case can be done by hand: if ${\cal P}$ is a maximal dps polytope with
size 2 in 
$\rr^2$, then  ${\cal L}({\cal P})$ must consist of four points chosen from
$$
\{(0,0), (0,1), (0,2), (1,0), (1,1), (2,0) \}.
$$ 
Since each congruence class is represented in  ${\cal L}({\cal P})$,
it must contain $(0,1)$, $(1,0)$ and $(1,1)$. These three points form a
parallelogram with each of the points $(0,0),\;(0,2)\,$ and $\,(2,0)$.  
Hence no fourth point can exist in ${\cal L}({\cal P})$
while preserving the dps property. The second case is similar, but
much more complicated.  Computer-aided calculations can be used to
conclude that no dps polytope in $\rr^3$ has size 4 or less. (We thank
Dr. Bruce Carpenter for doing the Mathematica coding.)

It can also be shown, using the style of argument 
of \cite[Ch.\ 3] {R1}, that every maximal dps polytope in $\rr^2$
is the image of the triangle in Example 1 under an affine unimodular
linear mapping, 
and consists of a triangle with area 3/2, and a single lattice point
inside, which will always be the centroid of the triangle. The
tetrahedron determined by ${\cal B}$ in Example 2 lies 
within the tetrahedron determined by ${\cal A}$, whereas in Example 3,
each point in $\cal L$ is on the boundary of the polytope. Thus there
are at least two distinct combinatorial types of maximal dps polytopes
in $\rr^3$.

We make no serious conjecture about the growth of $s_n$. On the one hand, any
maximal dps polytope must contain a lattice point with odd
coordinates, so $s_n \ge n$. In the other direction, it is not
difficult to use the proof of Theorem 3 to obtain a doubly-exponential
bound for $s_n$. Since this bound is likely to be very crude, we do
not present it explicitly. Another open question is to determine the minimum
volume of a maximal dps polytope in $\,{\mathbb R}^n\,$ for $\,n\geq
3$. We also do not know the answer to the following question: is every dps
polytope a subset of a maximal dps polytope?

We now discuss our original interest in this subject. Given $u
\in \zz_+^n$, define the monomial $x^u \in \rr[x_1,\dots,x_n]$ by 
$$
x^u = x_1^{u^{(1)}} \dots x_n^{u^{(n)}}.  
$$
Suppose ${\cal U} \sse \zz_+^n$ and consider the polynomial 
$$
p(x_1,\dots,x_{n}) = \sum_{u \in {\cal U}} b_u x^u.
$$
In \cite{CLR}, the present authors developed an algorithm for
determining whether $p$ 
can be written as a sum of squares of polynomials. A necessary
condition is that $p$ is psd; that is, $p(x_1,\dots,x_n) \ge 0$ for
all $x \in \rr^n$. Suppose $p$ is psd and let
$$
{\cal C}(p) = cvx\{ u: b_u \neq 0\}.
$$
Then ${\cal C}(p)$ is a lattice polytope; in fact it can be shown that  the
vertices of ${\cal C}(p)$ lie in $(2\zz)^n$, so that ${\cal P} := \frac 12
{\cal C}(p)$ is a lattice polytope. Let
$$
{\cal L}({\cal P}) = \{v_1,\dots, v_N\},
$$
and for $u \in {\cal C}(p)$, let $D(u) = \{(i,j): v_i + v_j = u\}$. 
It is proved in \cite[Thm.\ 2.4]{CLR} that $p$ can be written as a
sum of at most $r$ squares of polynomials if
and only if there is a real $N \times N$ symmetric psd matrix $A =
[a_{ij}]$ of rank at most $r$, so that 
$$
\sum_{(i,j) \in D(u) } a_{ij} = b_u\qquad \text{for all $u \in {\cal
C}(p)$.}
$$
If $\cal P$ is a dps polytope in $\rr^{n}$, then either $|D(u)| \le 1$
or $D(u) = \{(i,j),(j,i)\}$. In either case, $a_{ij}$ is completely
determined by $b_u$. In particular, if 
\begin{equation} \label{diag}
h_{{\cal P}}(x_1,\dots,x_{n}): = \sum_{i=1}^{N} \left(x^{v_i}\right)^2,
\end{equation}
then $A$ must equal $I_N$, the $N \times N$ identity matrix, so that $p$ is a
sum of $N$ squares, and no fewer.

Finally, we note that the homogenization of polynomials with $n$ variables
into forms with $n+1$ variables is precisely analogous to the embedding
of polytopes in $\,{\mathbb R}^n_+\,$ into the hyperplane $\Delta_{n+1}(s)$.

\medskip 

\noindent \textbf{Example 4.} (See \cite[Ex.\ 3.9]{CLR})

We return to Example 1, in its homogeneous version.
 Let $A = [a_{ij}]$ be a real symmetric $4
\times 4$ matrix and let
$$
f(t_1,t_2,t_3,t_4) = \sum_{i=1}^4  \sum_{j=1}^4 a_{ij}t_it_j
$$
be its associated quadratic form. We use the substitution suggested by 
${\cal L}({\cal P})$ and define the ternary sextic form
$$
p(x_1,x_2,x_3) = f(x_2x_3^2,x_1x_2^2,x_1^2x_3,x_1x_2x_3).
$$
Then $p$ is a sum of squares of polynomials (cubic forms) if and only
if $f$ is a psd quadratic form; that is,
$$
f(t_1,t_2,t_3,t_4) \ge 0\qquad \text{for all $(t_1,t_2,t_3,t_4) \in \rr^4$}.
$$
Since $t_4^3 = t_1t_2t_3$, the condition for $p$ to be a psd form is
weaker:
$$
f(t_1,t_2,t_3,(t_1t_2t_3)^{1/3}) \ge 0\qquad \text{for all
$(t_1,t_2,t_3) \in \rr^3$}. 
$$
If $f(t_1,t_2,t_3,t_4) = t_1^2
+ t_2^2 + t_3^2 - 3t_4^2$, then $f$ is not psd, but
$f(t_1,t_2,t_3,(t_1t_2t_3)^{1/3})\ge 0$ by the arithmetic-geometric
inequality. It follows that
$$
p(x_1,x_2,x_3) = x_2^2x_3^4 + x_1^2x_2^4 + x_1^4x_3^2 - 3x_1^2x_2^2x_3^2
$$
is a form which is psd, but not a sum of squares of polynomials.
This particular example was discussed in \cite{CL}. For a history and
bibliography of this subject and its relation to Hilbert's 17th
Problem, see  \cite{R2}.
 
\medskip
More generally, the {\it Pythagoras
number} of a ring $A$, $P(A)$,  is the smallest number $n \le \infty$
such that any sum of squares in $A$ can be expressed as a sum of at
most  $n$ squares in $A$.  Pfister \cite{P} proved in 1967 that
$P(\rr(x_1,\dots,x_n)) \le 2^n$. It is easy to see that  $P(\rr[x_1])
= 2$. Since maximal dps polytopes exist in $\rr^n$ for every $n$, a
consideration of  $h_{{\cal P}}$ (c.f.  \eqref{diag}) shows that
$P(\rr[x_1,\dots,x_n]) \ge 2^n$.  This is not the strongest result
possible: in \cite[p.60]{CDLR}, using other methods, Dai and the
present authors have shown that $P(\rr[x_1,\dots,x_n]) = \infty$ 
for $n \ge 2$.

\end{document}